\documentclass{birkjour}
 \usepackage{mathabx}
 \usepackage{xcolor}
 \newtheorem{thm}{Theorem}[section]
 
 \newtheorem{lem}[thm]{Lemma}
 \newtheorem{prop}[thm]{Proposition}
 \theoremstyle{definition}
 \newtheorem{defn}[thm]{Definition}
 \theoremstyle{remark}
 \newtheorem{rem}[thm]{Remark}
 
 \numberwithin{equation}{section}

\begin{document}

%-------------------------------------------------------------------------
% editorial commands: to be inserted by the editorial office
%
%\firstpage{1} \volume{228} \Copyrightyear{2004} \DOI{003-0001}
%\seriesextraline{This is the Concrete Title of this Book\br H.E. R and S.T.C. W, Eds.}
%\firstpage{1}

\title[Hankel Operators Between Köthe Spaces]
 {Hankel Operators Between Köthe Spaces}

%----------Author 1
\author{Nazl\i\;Do\u{g}an}
\address{Fatih Sultan Mehmet Vakif University, 
    34445 Istanbul, Turkey}
\email{ndogan@fsm.edu.tr ~\\
ORCID ID: 0000-0001-7548-954X}

%\thanks{This work was completed with the support of our
%----------Author 2
%\author{A Second Author}
%\address{The address of\br
%the second author\br
%sitting somewhere\br
%in the world}
%\email{dont@know.who.knows}
\subjclass{46A45, 47B35}

\keywords{Hankel Matrix, K\"othe Spaces, Compact Operators}

%\date{January 1, 2004}

\begin{abstract}
This paper is about the operators defined between K\"{o}the spaces whose associated matrix is a Hankel matrix. After demonstrating how these operators are defined, the conditions for continuity and compactness of these operators are constructed. It is shown that the backward and forward shift operators are mean ergodic and Ces\'{a}ro bounded by establishing a relationship between the backward and forward shift operators and Hankel and Toeplitz operators on power series spaces. 
\end{abstract}

%%% ----------------------------------------------------------------------
\maketitle
%%% ----------------------------------------------------------------------
%%%%%%%%%%%%%%%%%%%%%%%%%%%%%%%%%%%%%%%%%%%%%%%%%%%%%%%%%%%%%%%%%%%%%%%%%
\section{Introduction}

A finite or infinite matrix is called \emph{Hankel matrix} if its entries are constant along each skew-diagonal, that is, the matrix $(a_{m,n})$ is Hankel if $a_{m_{1},n_{1}} = a_{m_{2},n_{2}}$ whenever $m_{1} + n_{1} = m_{2} + n_{2}$. Hankel matrices are significant for several reasons and there are many applications in various fields such as; signal processing, control theory, and numerical analysis. Hankel matrices also serve as functional mathematical tools with diverse applications in engineering and computer science.

Hankel and Toeplitz operators defined on Hardy space of unit disk, $H^{2}(\mathbb{D})$, can be viewed as operators with infinite Hankel and infinite Toeplitz matrices, respectively, with respect to the standard orthonormal basis of 
 $H^{2}(\mathbb{D})$. A few years ago, Toeplitz operators, whose "associated" matrix is Toeplitz, are defined for more general topological vector spaces. For instance, in \cite{DJ}, Doma\'nski and Jasiczak developed the analogous theory for the space of $\mathcal{A}(\mathbb{R})$ real analytic functions on the real line.  In \cite{J}, Jasiczak introduced and characterized the class of Toeplitz operators
on the Fr\'echet space of all entire functions $\mathcal{O}(\mathbb{C})$.

In \cite{J}, Jasiczak defined  a continuous linear operator on $\mathcal{O}(\mathbb{C})$ as a Toeplitz operator if its matrix
is a Toeplitz matrix. The matrix of an operator is defined with
respect to the Schauder basis $(z^{n})_{n\in{N_{0}}}$.  The space of entire functions $\mathcal{O}(\mathbb{C})$ is isomorphic to a power series space of infinite type $\Lambda_{\infty}(n)$. By taking inspiration from Jasiczak paper \cite{J},  the author defined Toeplitz operators on more general power series spaces of finite or infinite type. In this paper, with the same idea, we will define the operators whose associated matrix is an infinite Hankel matrix between K\"{o}the spaces, especially power series spaces. We will construct some conditions for the continuity and compactness of these operators. In the final section, it is shown that the backward and forward shift operators are mean ergodic and Ces\'{a}ro bounded by establishing a relationship between the backward and forward shift operators and Hankel and Toeplitz operators on power series spaces. 

%%%%%%%%%%%%%%%%%%%%%%%%%%%%%%%%%%%%%%%%%%%%%%%%%%%%%%%%%%%%%%%%%%%%%%%%%
%%%%%%%%%%%%%%%%%%%%%%%%%%%%%%%%%%%%%%%%%%%%%%%%%%%%%%%%%%%%%%%%%%%%%%%%%
\section{Preliminaries}
In this section, after establishing terminology and notation, we collect some basic facts and definitions that are needed in the sequel. We will use the standard terminology and notation of \cite{MV}.

A complete Hausdorff locally convex space E whose topology is defined by a countable fundamental
system of seminorms $(\|\cdot\|_{k})_{k\in \mathbb{N}}$ is called a Fr\'echet space. A matrix $(a_{n,k})_{k,n\in \mathbb{N}}$ of non-negative numbers is called a Köthe matrix if it satisfies the following conditions:
\begin{itemize}
	\item[1.] For each $n\in \mathbb{N}$ there exists a $k\in \mathbb{N}$ with $a_{n,k}>0$. 
	\item[2.] $a_{n,k}\leq a_{n,k+1}$ for each $n,k\in \mathbb{N}$.
\end{itemize}
For a Köthe matrix $(a_{n,k})_{n,k\in \mathbb{N}}$,
\[ K(a_{n,k})=\bigg\{ x=(x_{n})_{n\in \mathbb{N}}: \;\|x\|_{k}:=\sum^{\infty}_{n=1}|x_{n}|a_{n,k}< \infty \quad\text{for all}\quad k\in \mathbb{N} \bigg\}\]
is called a Köthe space. Every Köthe space is a Fr\'echet space. From Proposition 27.3 of \cite{MV}, the dual space of a Köthe space is
\[ (K(a_{n,k}))^{\prime}=\bigg\{ y=(y_{n})_{n\in \mathbb{N}}\; \bigg|\; \sup_{n\in \mathbb{N}}|y_{n}a_{n,k}^{-1}|<+\infty \; \text{for some}  \; k\in \mathbb{N}\bigg\}.\]

Grothendieck-Pietsch Criteria \cite{MV}[Theorem 28.15] states that a Köthe space $K(a_{n,k})$ is nuclear if and only if  for every $k\in \mathbb{N}$, there exists a $l>k$ so that $$\displaystyle \sum^{\infty}_{n=1}{a_{n,k}\over a_{n,l}}<\infty.$$
For a nuclear Köthe space, $\|x\|_{k}=\sup_{n\in \mathbb{N}}|x_{n}|a_{n,k}$, $k\in \mathbb{N}$ forms an equivalent system of seminorms to the fundamental system of seminorms $\|x\|_{k}=\sum^{\infty}_{n=1}|x_{n}|a_{n,k}$, $k\in \mathbb{N}$.

Let $\alpha=\left(\alpha_{n}\right)_{n\in \mathbb{N}}$ be a non-negative increasing sequence with $\displaystyle \lim_{n\rightarrow \infty} \alpha_{n}=\infty$. A power series space of finite type is defined by
$$\Lambda_{1}\left(\alpha\right):=\left\{x=\left(x_{n}\right)_{n\in \mathbb{N}}: \;\left\|x\right\|_{k}:=\sum^{\infty}_{n=1}\left|x_{n}\right|e^{-{1\over k}\alpha_{n}}<\infty \textnormal{ for all } k\in \mathbb{N}\right\}$$
and a power series space of infinite type is defined by
$$\displaystyle \Lambda_{\infty}\left(\alpha\right):=\left\{x=\left(x_{n}\right)_{n\in \mathbb{N}}:\; \left\|x\right\|_{k}:=\sum^{\infty}_{n=1}\left|x_{n}\right|e^{k\alpha_{n}}<\infty \textnormal{ for all } k\in \mathbb{N}\right\}.$$
Power series spaces form an important family of Köthe spaces  and they contain the spaces of holomorphic functions on $\mathbb{C}^{d}$ and $\mathbb{D}^{d}$,
\[\mathcal{O}(\mathbb{C}^{d})\cong \Lambda_{\infty}(n^{\frac{1}{d}}) \quad\quad\text{and} \quad\quad \mathcal{O}(\mathbb{D}^{d})\cong \Lambda_{1}(n^{\frac{1}{d}})
\]
where $\mathbb{D}$ is the unit disk in $\mathbb{C}$ and $d\in \mathbb{N}$.
%%%%%%%%%%%%%%%%%%%%%%%%%%%%%%%%%%%%%%%%%%%%%%%%%%%%%%%%%%%%%%%%%%%%%%%%

Let $E$ and $F$ be Fr\'echet spaces. A linear map $T:E\to F$ is called continuous 
if for every $k\in \mathbb{N}$ there exists $p\in \mathbb{N}$ 
and $C_k>0$ such that 
\[\|Tx\|_k\leq C_k\|x\|_p\]
for all $x\in E$. A linear map $T: E\to F$ is called compact
if $T(U)$ is precompact in $F$ where $U$ is a neighborhood of zero of E. 

In this paper, we fixed the symbol $e_{n}$ to denote the sequence \[(0,0,\dots,0,1,0,\dots)\] where 1 is in the n$^{th}$ place and 0 is in the others.

We will use the following Lemma to determine the continuity and compactness of operators defined between Köthe spaces. 
%%%%%%%%%%%%%%%%%%%%%%%%%%%%%%%%%%%%%%%%%%%%%%%%%%%%%%%%%%%%%%%%%%%%%%%%%
\begin{lem}\label{Crone} 
	Let $K(a_{n,k})$ and $K(b_{n,k})$ be Köthe spaces.
	\begin{itemize}
		\item[a.] $T: K(a_{n,k}) \to K(b_{n,k})$ is a linear continuous operator if and only if for each k there exists m such that
		\[ \sup_{n\in \mathbb{N}} {\frac{\|Te_{n}\|_{k}}{\|e_{n}\|_{m}}<\infty }.\]
		\item[b.] If $K(b_{n,k})$ is Montel, then $T: K(a_{n,k}) \to K(b_{n,k})$ is a compact operator if and only if there exists m such that for all k
		\[\sup_{n\in \mathbb{N}} {\frac{\|Te_{n}\|_{k}}{\|e_{n}\|_{m}}<\infty }.\]
	\end{itemize}
\end{lem}
\begin{proof} Lemma 2.1 of \cite{CR75}.
\end{proof}
%%%%%%%%%%%%%%%%%%%%%%%%%%%%%%%%%%%%%%%%%%%%%%%%%%%%%%%%%%%%%%%%%%%%%%%%%%%%%%%%%%%%%%%%%%%%
A Fr\'echet space $E$ is Montel if each bounded set in $E$ is relatively compact. Every power series space is Montel, see Theorem 27.9 of \cite{MV}.
%%%%%%%%%%%%%%%%%%%%%%%%%%%%%%%%%%%%%%%%%%%%%%%%%%%%%%%%%%%%%%%%%%%%%%%%

The next proposition says that the continuity condition is sufficient to ensure that linear operators defined only on the basis elements are well-defined.
\begin{prop}\label{T10} Let $K(a_{n,k})$, $K(b_{n,k})$ be Köthe spaces and  $(a_{n})_{n\in \mathbb{N}}\in K(b_{n,k})$  be a sequence. Let us define a linear map $T:K(a_{n,k})\to K(b_{n,k})$ such as
	\[Te_{n}=a_{n} \hspace{0.5in}\text{and}\hspace{0.5in} Tx=\sum^{\infty}_{n=1}x_{n}Te_{n}\]
	for every $x=\sum^{\infty}_{n=1} x_{n}e_{n}$ and $n\in \mathbb{N}$. If the  continuity condition 
	\[\forall k\in \mathbb{N} \quad\quad \exists m\in \mathbb{N}\quad\quad\quad\quad \sup_{n\in \mathbb{N}}\frac{ \|Te_{n}\|_{k}}{\|e_{n}\|_{m}}<\infty \]
	holds, then $T$ is well-defined and continuous operator.
\end{prop}
\begin{proof} Proposition 2.2 of \cite{N}.
\end{proof}

In this paper, we will call an operator which is defined between Köthe spaces as a Hankel operator if its matrix is a Hankel matrix defined with respect to the Schauder basis $(e_{n})_{n\in \mathbb{N}}$. We will concentrate on Hankel operators defined between power series spaces and determine the conditions that give us the continuity and compactness of these operators.
%%%%%%%%%%%%%%%%%%%%%%%%%%%%%%%%%%%%%%%%%%%%%%%%%%%%%%%%%%%%%%%%%%%%%%%%%
\section{Hankel Operator Defined Between K\"{o}the Spaces}
%%%%%%%%%%%%%%%%%%%%%%%%%%%%%%%%%%%%%%%%%%%%%%%%%%%%%%%%%%%%%%%%%%%%%%%%%
Let $\theta=(\theta_{n})_{n\in \mathbb{N}_{0}}$ be any sequence. The Hankel matrix defined by $\theta$ is
\[
H_{\theta}=\begin{pmatrix}
\theta_{0} &\theta_{1}&\theta_{2}&\theta_{3}&\cdots \\
\theta_{1}&\theta_{2}&\theta_{3}&\theta_{4}&\cdots\\
\theta_{2}&\theta_{3}&\theta_{4}&\theta_{5}&\cdots\\
\theta_{3}&\theta_{4}&\theta_{5}&\theta_{6}&\cdots\\
\vdots&\vdots&\vdots&\vdots&\ddots
\end{pmatrix}.
\]
We aim to define an operator $H_{\theta}: K(a_{n,k})\to K(b_{n,k})$ by taking  $H_{\theta}e_{n}$ as the n$^{th}$ column of the above matrix, that is,
\[H_{\theta}e_{n}=(\theta_{n-1},\theta_{n},\theta_{n+1},\cdots)=\sum^{\infty}_{j=1}\theta_{j+n-2}e_{j}\]
provided that $H_{\theta}e_{n}\in K(b_{n,k})$ for every $n\in \mathbb{N}$. 
Therefore, for every $x=\sum^{\infty}_{n=1}x_{n}e_{n}\in K(a_{n,k})$, the operator  $H_{\theta}$ can be written as 
\begin{equation}\label{E1} H_{\theta}x=\sum^{\infty}_{n=1}x_{n} H_{\theta}e_{n}.\end{equation}
Actually, we cannot confirm that the operator $H_{\theta}:K(a_{n,k})\to K(b_{n,k})$ is well-defined as we cannot guarantee that the series $\sum^{\infty}_{n=1}x_{n} \widehat{T}_{\theta}e_{n}$ converges in $K(b_{n,k})$ for every $x\in K(a_{n,k})$. In this section, we will share some conditions under which this operator can be appropriately defined between power series spaces and in those instances, we will analyze its continuity and compactness.

%------------------------------------------------------------------------------------------------------------------------------------

As a direct consequence of Proposition \ref{T10}, we have the following:
\begin{prop}\label{P0} $H_{\theta}:K(a_{n,k})\to K(b_{n,k})$ is well-defined and continuous if and only if $H_{\theta}e_{n}\in K(b_{n,k})$ for every $n\in \mathbb{N}$ and the continuity condition 
\[\forall k\in \mathbb{N} \quad\quad \exists m\in \mathbb{N} \quad\quad\quad\quad \sup_{n\in \mathbb{N}}\frac{ \|H_{\theta}e_{n}\|_{k}}{\|e_{n}\|_{m}}<\infty \]
holds.
\end{prop}
%----------------------------------------------------------------------------------------------------------------
\begin{rem}\label{l1} 
When $H_{\theta}$ defines a continuous linear operator, it can be especially said that 
\[H_{\theta}e_{1}=(\theta_{0},\theta_{1}, \theta_{2},\cdots)=\theta\]
is always in $K(b_{n,k})$.
\end{rem}
%%%%%%%%%%%%%%%%%%%%%%%%%%%%%%%%%%%%%%%%%%%%%%%%%%%%%%%%%%%%%%%%%%%%%%%%%%%%%%%%%%%%%%%%%%%%%%%%%%%%%%%%%%%%%%

As mentioned in  Remark \ref{l1} the sequence  $\theta$ lies in the range space of Hankel operator $H_{\theta}$. In the proposition below we will demonstrate that the sequence $\theta$ should be in the dual space of the domain space of Hankel operator $H_{\theta}$.

\begin{prop}\label{P1}
Let $K(a_{n,k})$, $K(b_{n,k})$ be Köthe spaces. If  $H_{\theta}: K(a_{n,k})\to K(b_{n,k})$ is a continous linear operator, then 
$\theta\in (K(a_{n,k}))^{\prime}$.
\end{prop}
\begin{proof} Let $H_{\theta} :K(a_{n,k}) \to K(b_{n,k})$ be a continuous linear operator with the formula 
\[H_{\theta}e_{n}=(\theta_{n-1},\theta_{n},\theta_{n+1},\cdots)=\sum^{\infty}_{j=1}\theta_{j+n-2}e_{j}.\]
for every $n\in \mathbb{N}$. By Lemma \ref{Crone}, for all $k\in \mathbb{N}$ there exist $m\in \mathbb{N}$ and $C>0$ such that
$$\quad\quad\quad \quad \|H_{\theta}e_{n}\|_{k}=\sum^{\infty}_{j=1}|\theta_{j+n-2}|b_{j,k}\leq C\|e_{n}\|_{m}=Ca_{n,m}\quad\quad\quad \quad \quad\forall n\in \mathbb{N}.$$
Then, we have that for all $n, j\in \mathbb{N}$
\begin{equation}|\theta_{j+n-2}|b_{j,k}\leq Ca_{n,m}.
\end{equation}
Since $(b_{n,k})_{n,k\in \mathbb{N}}$ is Köthe matrix, there exists a $k_{0}\in \mathbb{N}$ such that $b_{1,k_{0}}\neq 0$. Hence there exist $m_{0}\in \mathbb{N}$ and $C_{0}>0$ such that
\[ \hspace{0.75in}\quad \quad\quad \quad \quad\quad|\theta_{n-1}|\leq \frac{C_{0}}{b_{1,k_{0}}}a_{n,m_{0}}\quad\quad\quad\quad \quad \quad\hspace{0.55in}\forall n\in \mathbb{N}.\]
This says that $\theta\in (K(a_{n,k}))^{\prime}$.
\end{proof}
%%%%%%%%%%%%%%%%%%%%%%%%%%%%%%%%%%%%%%%%%%%%%%%%%%%%%%%%%%%%%%%%%%%%%%%%%
%-----------------------------------------------------------------------------------------------------------------------------------

%--------------------------------------------------------------------------------------------------------------------------------------------------------
\subsection{Power Series Spaces as a Range Space}
%-----------------------------------------------------------------------------------------------------------------------------------------------------
In this subsection, we examine the continuity and compactness of Hankel operators from a Köthe space to a power series space. Initially, we will start with the case that the range space is a power series space of infinite type $\Lambda_{\infty}(\beta)$.
 
%%%%%%%%%%%%%%%%%%%%%%%%%%%%%%%%%%%%%%%%%%%%%%%%%%%%%%%%%%%%%%%%%%%%%%%%%%%%%%%%%%%%%%%%%%%
\begin{prop}\label{P2} Let $K(a_{n,k})$ be a Köthe space and $\theta\in \Lambda_{\infty}(\beta)$. If the condition
\begin{equation}\label{11}\hspace{0.7in}
\exists  m_{0}\in \mathbb{N}, \; C>0 \hspace{0.5in}a_{n,m_{0}}\geq C \hspace{0.9in} \forall n\in \mathbb{N}
\end{equation}
holds, then $H_{\theta}:K(a_{n,k}) \to \Lambda_{\infty}(\beta)$ is well-defined, continuous and compact.
\end{prop}
\begin{proof} Let us assume that $\theta\in \Lambda_{\infty}(\beta)$ and the condition (\ref{11}) holds.  This gives us that for every $k\in \mathbb{N}$ there exists a $D>0$ such that
\begin{equation*}
\begin{split}
\|H_{\theta}e_{n}\|_{k}& =\sum^{\infty}_{j=1}|\theta_{j+n-2}|e^{k\beta_{j}}\leq \sum^{\infty}_{j=1}|\theta_{j+n-2}|e^{k\beta_{j+n-1}} \\ &\leq \|\theta\|_{k}\leq \frac{1}{C} \|\theta\|_{k}a_{n,m_{0}}\leq D a_{n,m_{0}}
\end{split}
\end{equation*}
for every $n\in \mathbb{N}$. This says that $H_{\theta}e_{n}\in K(b_{n,k})$ for every $n\in \mathbb{N}$ and for every $k\in \mathbb{N}$
$$\sup_{n\in \mathbb{N}} \frac{\|H_{\theta}e_{n}\|_{k}}{\|e_{n}\|_{m_{0}}}<\infty.$$
Then $H_{\theta}$ is well-defined and continuous from Proposition \ref{P0}. Since every power series space is Montel and $m_{0}$ does not depend on k, $H_{\theta}$ is also compact from Lemma \ref{Crone}.
\end{proof}

%------------------------------------------------------------------------------------------------------------

By a direct consequence of Proposition \ref{P2} we can give the following theorem:

\begin{thm}\label{T1} For every $\theta\in \Lambda_{\infty}(\beta)$, the Hankel operator $H_{\theta}$ from any infinite type power series space $\Lambda_{\infty}(\alpha)$ to $\Lambda_{\infty}(\beta)$ is continous and compact. 
\end{thm}
%--------------------------------------------------------------------------------------------------
Now we want to write a weaker condition on the matrix of $K(a_{n,k})$ in Proposition \ref{P2}. %However to achieve this, we require the following growth condition
%on the sequence $\beta$:
%\begin{equation}\label{E2}\forall s,t\quad\quad\quad\quad \beta_{t}+\beta_{s}\leq \beta_{t+s-1}\end{equation}
%for $\Lambda_{\infty}(\beta).$
%-------------------------------------------------------------------------------------------------------------------------------------------
\begin{prop}\label{P3}  Let $K(a_{n,k})$ be a Köthe space and $\theta\in \Lambda_{\infty}(\beta)$.
% and $\beta$ are a nonnegative increasing sequence that tends to infinity and satisfies the condition in \ref{E2}.
Assume that  the following condition holds:
\begin{equation}\label{E3}\quad\forall k\in \mathbb{N}\quad\quad \exists m\in \mathbb{N}, C>0\quad \quad\quad\quad e^{-k\beta_{n}}\leq C a_{n,m}\quad\quad\quad \forall n\in \mathbb{N}. \end{equation}
Then $H_{\theta}:K(a_{n,k}) \to \Lambda_{\infty}(\beta)$ is well-defined and continuous.
\end{prop}

\begin{proof} Let $\theta\in \Lambda_{\infty}(\beta)$. Since $\beta$ is increasing, $\max\{\beta_{j},\beta_{n}\}\leq \beta_{j+n-1}$ and $\beta_{j}+\beta_{n}\leq 2\beta_{j+n-1}$ for all $j,n\in \mathbb{N}$.
By using the condition (\ref{E3}), for every $k\in \mathbb{N}$ there exist $m\in \mathbb{N}$ and $C>0$ such that
\begin{equation*}
\begin{split}
\|H_{\theta}e_{n}\|_{k}&=\sum^{\infty}_{j=1}|\theta_{j+n-2}|e^{k\beta_{j}}\leq \sum^{\infty}_{j=1}|\theta_{j+n-2}|e^{2k\beta_{j+n-1}}e^{k(\beta_{j}-2\beta_{j+n-1})}\\
&\leq \sum^{\infty}_{j=1}|\theta_{j+n-2}|e^{2k\beta_{j+n-1}}e^{-k\beta_{n}}=\|\theta\|_{k}e^{-k\beta_{n}}\leq C\|\theta\|_{k}a_{m,n}
\end{split}
\end{equation*}
and then $H_{\theta}e_{n}\in \Lambda_{\infty}(\beta)$ for every $n \in \mathbb{N}$ and 
\[\sup_{n\in \mathbb{N}}\frac{ \|H_{\theta}e_{n}\|_{k}}{\|e_{n}\|_{m}}<\infty.\]
Proposition \ref{P0} says that $H_{\theta}$ is well-defined and  continuous.
\end{proof}
%----------------------------------------------------------------------------------------------------------------------------------------------------------------------
\begin{prop}\label{P4} Let $K(a_{n,k})$ be a Köthe space and $\theta\in \Lambda_{\infty}(\beta)$.
Assume that the following condition holds:
\begin{equation}\label{E4}\quad\quad\exists m\in \mathbb{N} \quad\forall k\in \mathbb{N}\quad \exists C>0 \quad\quad\quad e^{-k\beta_{n}}\leq C a_{n,m}\quad\quad\quad \forall n\in \mathbb{N}.\end{equation}
Then $H_{\theta}:K(a_{n,k}) \to \Lambda_{\infty}(\beta)$ is compact.
\end{prop}
\begin{proof} Substituting condition \ref{E3} with condition \ref{E4} and following the steps outlined in Proposition \ref{P3}, it can be given as a consequence of Lemma \ref{Crone}.
\end{proof}
%----------------------------------------------------------------------------------------------------------

As a direct consequence of Proposition \ref{P3} and Proposition \ref{P4} we can give the following theorem:

\begin{thm}\label{T2} Let $\beta,\alpha$ be two nonnegative increasing sequences that tend to infinity. Assume that there exist $A,B>0$ such that
\begin{equation}\label{C1}
\alpha_{n}\leq A\beta_{n} +B
\end{equation}
for all $n\in \mathbb{N}$. Then, for every $\theta\in \Lambda_{\infty}(\beta)$, the Hankel operator $H_{\theta}: \Lambda_{1}(\alpha) \to \Lambda_{\infty}(\beta)$ is well-defined, continous and  compact. 
\end{thm}
\begin{proof} Let us assume that there exist $A,B>0$ satisfying $\alpha_{n}\leq A\beta_{n} +B$ for all $n\in \mathbb{N}$. Then for all $m,k\in \mathbb{N}$ we write
$$ \frac{1}{mA}\alpha_{n}-\frac{B}{A}\leq \frac{1}{A}\alpha_{n}-\frac{B}{A}\leq \beta_{n}\leq k\beta_{n}$$
and 
$$-k\beta_{n}\leq \frac{B}{A}-\frac{1}{mA}\alpha_{n}$$
for all $n\in \mathbb{N}$.
Then for all $\tilde{m}\in \mathbb{N}$ satisfying $\tilde{m}>mA$ and for all $k\in \mathbb{N}$, there exists a $C>0$ such that
$$e^{-k\beta_{n}}\leq C e^{-\frac{1}{\tilde{m}}\alpha_{n}}$$
for every $n\in \mathbb{N}$. This says that the conditions in (\ref{E3}) and (\ref{E4}) are satisfied. From Proposition \ref{P3} and \ref{P4}, $H_{\theta}:\Lambda_{1}(\alpha)\to \Lambda_{\infty}(\beta)$ is well-defined, continuous and compact.
\end{proof}
%-----------------------------------------------------------------------------------------------------

%--------------------------------------------------------------------------------------------------
Now, we will explore the continuity and compactness of the Hankel operator $H_{\theta}$, when the range space is a power series space of finite type $\Lambda_{1}(\beta)$. To this, we require the stability condition on the sequence $\beta$. A sequence $\beta$ is called stable if
\begin{equation}\label{E3.5}
	\sup_{n\in \mathbb{N}} {\beta_{2n}\over \beta_{n}}<\infty.
\end{equation}
%we require the following growth condition
%on the sequence $\beta$:
%\begin{equation}\label{E5}
%\forall s,t\quad\quad\quad\quad \beta_{t+s-1}\leq \beta_{t}+\beta_{s}.
%\end{equation}

%--------------------------------------------------------------------------------------
\begin{prop}\label{P5}
Let $\beta$ be a stable sequence,  $K(a_{n,k})$ be a Köthe space and  $\theta\in \Lambda_{1}(\beta)$.
Assume that the following condition holds:
\begin{equation}\label{E6}\quad \forall k\in \mathbb{N}\quad\quad \exists m\in \mathbb{N}, C>0 \quad\quad\quad e^{\frac{1}{k}\beta_{n}}\leq C a_{m,n}\quad\quad\quad\quad\quad \forall n\in \mathbb{N}. \end{equation}
Then $H_{\theta}:K(a_{n,k})\to \Lambda_{1}(\beta)$ is well-defined and continuous. 
\end{prop}

\begin{proof} 
Let $\theta\in \Lambda_{1}(\beta)$.  Since $\beta$ is stable, there exists an $M> 1$ such that
\[\hspace{1.75in}\beta_{2n}\leq M\beta_{n}\hspace{1.75in} \forall n\in \mathbb{N}.\]
Since $\beta$ is increasing we have the following: Let $j,n\in \mathbb{N}$. If $j+n-1=2t$ or $j+n-1=2t-1$ for some $t\in \mathbb{N}$, then
$$\beta_{j+n-1}\leq \beta_{2t}\leq M\beta_{t}\leq M\max\{ \beta_{n},\beta_{j}\} \leq M(\beta_{n}+\beta_{j})$$
since $t\leq \max\{n,j\}$ and $\beta$ is increasing. Therefore, we have 
\begin{equation}\label{s1}
\beta_{j+n-1}\leq M(\beta_{n}+\beta_{j})
\end{equation}
for all $j,n\in \mathbb{N}$. By using the condition in \ref{E6}, we can write that for every $k\in \mathbb{N}$ there exist $m\in \mathbb{N}$ and $C>0$ such that
\begin{equation*}
\begin{split}
\|H_{\theta}e_{n}\|_{k}&=\sum^{\infty}_{j=1}|\theta_{j+n-2}|e^{-\frac{1}{k}\beta_{j}}\leq \sum^{\infty}_{j=1}|\theta_{j+n-2}|e^{-\frac{1}{Mk}\beta_{j+n-1}}e^{\frac{1}{k}(\frac{1}{M}\beta_{j+n-1}-\beta_{j})}\\
&\leq e^{\frac{1}{k}\beta_{n}}\sum^{\infty}_{j=1}|\theta_{j+n-2}|e^{-\frac{1}{Mk}\beta_{j+n-1}}\leq e^{\frac{1}{k}\beta_{n}} \|\theta\|_{Mk}\leq C a_{m,n}
\end{split}
\end{equation*}
and then $H_{\theta}e_{n}\in \Lambda_{1}(\beta)$ for every $n \in \mathbb{N}$ and 
\[\sup_{n\in \mathbb{N}}\frac{ \|H_{\theta}e_{n}\|_{k}}{\|e_{n}\|_{m}}<\infty.\]
From Proposition \ref{P0}, $H_{\theta}$ is well-defined and continuous.
\end{proof}
%-----------------------------------------------------------------------------------------------
 By modifying the condition in (\ref{E6}), we can establish a condition for the compactness of the operator $H_{\theta}: K(a_{n,k}) \to \Lambda_{1}(\beta)$.
%------------------------------------------------------------------------------------------------------------------------------------------------------------
\begin{prop}\label{P6}
Let $\beta$ be a stable sequence, $K(a_{n,k})$ be a Köthe space and $\theta\in \Lambda_{1}(\beta)$.
Assume that the following condition holds:
\begin{equation}\label{E7}\exists m\in \mathbb{N} \quad\forall k\in \mathbb{N}\quad\exists C>0 \quad\quad\quad e^{\frac{1}{k}\beta_{n}}\leq C a_{m,n}\quad\quad\quad\quad \forall n\in \mathbb{N}. \end{equation}
Then $H_{\theta}:K(a_{n,k})\to \Lambda_{1}(\beta)$ is compact.
\end{prop}
\begin{proof} It follows similar steps to the proof of Proposition \ref{P5}.
\end{proof}
%------------------------------------------------------------------------------------------------------------------------------------------------------------------------------------------------------------ 
As a consequence of Proposition \ref{P3} and Proposition \ref{P4} we can give the following theorem:

\begin{thm}\label{T3} Let $\beta,\alpha$ be two nonnegative increasing sequences that tend to infinity. Assume that $\beta$ is stable and there exist $A,B>0$ such that
\begin{equation}\label{C2}\beta_{n}\leq A\alpha_{n} +B
\end{equation}
for all $n\in \mathbb{N}$. Then, for every $\theta\in \Lambda_{1}(\beta)$, the Hankel operator $H_{\theta}: \Lambda_{\infty}(\alpha) \to \Lambda_{1}(\beta)$ is well-defined, continous and compact.
\end{thm}
\begin{proof} Let assume that there exist $A,B>0$ satisfying $\beta_{n}\leq A\alpha_{n} +B$ for all $n\in \mathbb{N}$. Then for all $m,k\in \mathbb{N}$ we write	
$$ \frac{1}{k}\beta_{n}\leq \beta_{n}\leq A\alpha_{n}+B \leq mA\alpha_{n}+B$$
for all $n\in \mathbb{N}$. Then for all $\tilde{m}\in \mathbb{N}$ satisfying $\tilde{m}>mA$ and for all $k\in \mathbb{N}$, there exists a $C>0$ such that
$$e^{\frac{1}{k}\beta_{n}}\leq C e^{\tilde{m}\alpha_{n}}$$
for every $n\in \mathbb{N}$. This says that the conditions in (\ref{E6}) and (\ref{E7}) are satisfied. From Proposition \ref{P5} and \ref{P6}, $H_{\theta}:\Lambda_{1}(\alpha)\to \Lambda_{\infty}(\beta)$ is well-defined, continuous and compact.
\end{proof}
%-------------------------------------------------------------------------------------------------------------------------------------------------------
\subsection{Power Series Spaces as a Domain Space}
%----------------------------------------------------------------------------------------------------------------------------------------------------- 

In this subsection, we examine the continuity and compactness of Hankel operators from a power series space to a K\"{o}the space. In the first step, we will consider the operators $H_{\theta}:\Lambda_{\infty}(\alpha)\to K(b_{n,k})$ for $\theta\in (\Lambda_{\infty}(\alpha))^{\prime}$.

\begin{prop}\label{P7}
 Let $\alpha$ be a stable sequence, $K(b_{n,k})$ be a Köthe space and $\theta\in (\Lambda_{\infty}(\alpha))^{\prime}$. Assume that the following condition holds: 
\begin{equation}\label{E8}\forall m\in \mathbb{N}\quad\quad e^{m\alpha_{n}}\in K(b_{n,k}). 
\end{equation}
In this case $H_{\theta}:\Lambda_{\infty}(\alpha)\to K(b_{n,k})$ is well-defined and continuous. If $K(b_{n,k})$ is Montel, then $H_{\theta}$ is also compact.
\end{prop}
\begin{proof} Let $\theta\in (\Lambda_{\infty}(\alpha))^{\prime}$. Then there exist $m_{0}\in \mathbb{N}$ and $C>0$  such that
\[|\theta_{n-1}|\leq Ce^{m_{0}\alpha_{n}}\]
for every $n\in \mathbb{N}$. Since $\alpha$ is stable, there exists a $M>0$ such that
$$\alpha_{j+n-1}\leq M(\alpha_{n}+\alpha_{j})$$
for all $j,n\in \mathbb{N}$ as we proved in (\ref{s1}) of the proof of Proposition \ref{P5}. By using the condition in (\ref{E8}), we can write that for every $k\in \mathbb{N}$ there exists a $\widetilde{C}>0$ such that
\begin{equation*}
\begin{split}
\|H_{\theta}e_{n}\|_{k}&=\sum^{\infty}_{j=1}|\theta_{j+n-2}|b_{n,k}\leq C \sum^{\infty}_{j=1}e^{m_{0}\alpha_{j+n-1}} b_{k,j} \\
&=Ce^{m_{0}M\alpha_{n}} \sum^{\infty}_{j=1} e^{m_{0}M\alpha_{j}}b_{k,j}= C\|e_{n}\|_{m_{0}M}\|e^{m_{0}M\alpha_{n}}\|_{k} =\widetilde{C}\|e_{n}\|_{m_{0}M}
\end{split}
\end{equation*}
this means that $H_{\theta} e_{n}\in K(b_{n,k})$ for every $n\in \mathbb{N}$ and for all $k\in \mathbb{N}$ 
\[\sup_{n\in \mathbb{N}}\frac{ \|H_{\theta}e_{n}\|_{k}}{\|e_{n}\|_{m_{0}M}}<\infty.\]
Proposition \ref{T10} says that $H_{\theta}$ is well defined and continuous. Since $m_{0}M$ does not depend on k, $H_{\theta}$ is compact provided that $K(b_{n,k})$ is Montel. 
\end{proof}
%--------------------------------------------------------------------
As a consequence of Proposition \ref{P7}, we have the following:
\begin{thm} Let $\alpha$ be a stable sequence and $\Lambda_{1}(\beta)$ be a nuclear power series space of finite type. If the following condition
\[n\alpha_{n}\leq \beta_{n}\]
holds for all sufficiently large $n\in \mathbb{N}$. Then, for every $\theta\in (\Lambda_{\infty}(\alpha))^{\prime}$, the Hankel operator $H_{\theta}: \Lambda_{\infty}(\alpha) \to \Lambda_{1}(\beta)$ is well-defined, continous and compact.  
\end{thm}
\begin{proof} For every $m,k\in \mathbb{N}$ and for all sufficiently large $n$, we have
\[mk\alpha_{n}\leq n\alpha_{n}\leq \beta_{n}.\]
 This gives us that
\[\sup_{n\in \mathbb{N}} e^{m\alpha_{n}-{1\over k}\beta_{n}}<+\infty.\]
This means that $e^{m\alpha_{n}}\in \Lambda_{1}(\beta)$ for every $m\in \mathbb{N}$. Proposition \ref{P8} says that $H_{\theta}:\Lambda_{\infty}(\alpha)\to \Lambda_{1}(\beta)$ is well-defined, continuous and compact for every $\theta\in (\Lambda_{\infty}(\alpha))^{\prime}$.
\end{proof}
%------------------------------------------------------------------------------------------------------------------------------------------------------

Now we will consider the operators $H_{\theta}:\Lambda_{1}(\alpha)\to K(b_{n,k})$ for $\theta\in (\Lambda_{1}(\alpha))^{\prime}$.

\begin{prop}\label{P8} Let $\theta\in (\Lambda_{1}(\alpha))^{\prime}$.
Assume that the following condition holds: 
\begin{equation}\label{E9}\forall m\in \mathbb{N}\quad\quad e^{-\frac{1}{m}\alpha_{n}}\in K(b_{n,k}). \end{equation}
Then $H_{\theta}:\Lambda_{1}(\alpha)\to K(b_{n,k})$ is well-defined and continuous. If $K(b_{n,k})$ is Montel, then $H_{\theta}$ is compact.
\end{prop}
\begin{proof} Since $\alpha$ is increasing, $\max\{\alpha_{j},\alpha_{n}\}\leq \alpha_{j+n-1}$ and $\alpha_{j}+\alpha_{n}\leq 2\alpha_{j+n-1}$ for all $j,n\in \mathbb{N}$.	 Let $\theta\in (\Lambda_{1}(\alpha))^{\prime}$. Then there exist $m_{0}\in \mathbb{N}$ and $C>0$  such that
\[|\theta_{n-1}|\leq Ce^{-\frac{1}{m_{0}}\alpha_{n}}\]
for every $n\in \mathbb{N}$. 
By using this and the condition in (\ref{E8}), we can write that for every $k\in \mathbb{N}$ there exists a $\widetilde{C}>0$ such that
\begin{equation*}
\begin{split}
\|H_{\theta}e_{n}\|_{k}&=\sum^{\infty}_{j=1}|\theta_{j+n-2}|b_{n,k}\leq C \sum^{\infty}_{j=1}e^{-\frac{1}{m_{0}}\alpha_{j+n-1}} b_{k,j} \\
&=Ce^{-\frac{1}{2m_{0}}\alpha_{n}} \sum^{\infty}_{j=1} e^{-\frac{1}{2m_{0}}\alpha_{j}}b_{k,j}= C\|e_{n}\|_{2m_{0}}\|e^{-\frac{1}{2m_{0}}\alpha_{n}}\|_{k} =\widetilde{C}\|e_{n}\|_{2m_{0}}
\end{split}
\end{equation*}
and then
\[\sup_{n\in \mathbb{N}}\frac{ \|H_{\theta}e_{n}\|_{k}}{\|e_{n}\|_{2m_{0}}}<\infty.\]
$H_{\theta}$ is well-defined and continuous by Proposition \ref{T10}. Since $m_{0}$ does not depend on $k$, $H_{\theta}$ is compact provided that $K(b_{n,k})$ is Montel. 
\end{proof}
Proposition \ref{P8} gives us the following result:
%---------------------------------------------------------------------------------------------------------------------------------------------------
\begin{thm}\label{T4} Let $\Lambda_{1}(\beta)$ be a nuclear power series spaces of finite type. For every $\theta\in (\Lambda_{1}(\alpha))^{\prime}$, the Hankel operator $H_{\theta}: \Lambda_{1}(\alpha) \to \Lambda_{1}(\beta)$ is well-defined, continous and compact.
\end{thm}
\begin{proof} Since 
$$\sup_{n\in \mathbb{N}} e^{-\frac{1}{m}\alpha_{n}}e^{-\frac{1}{k}\beta_{n}}<+\infty$$
for every $m,k\in \mathbb{N}$, then $e^{-\frac{1}{m}\alpha_{n}}\in \Lambda_{1}(\beta)$ for every $m\in \mathbb{N}$. Proposition \ref{P8} says that $H_{\theta}:\Lambda_{1}(\alpha)\to \Lambda_{1}(\beta)$ is well-defined, continuous and compact for every  $\theta\in (\Lambda_{1}(\alpha))^{\prime}$.
\end{proof}

%---------------------------------------------------------------------------------------------------------------------------------------------------------------------------------------------------------
\subsection{S-Tameness of The Family of Hankel Operators}
A grading on a Fr\'echet space $E$ consists of a sequence of seminorms $\{\|\cdot\|_{n}\}_{n\in \mathbb{N}}$ that are increasing, which means that for every $x\in E$, the inequalities 
$\|x\|_{1}\leq \|x\|_{2}\leq \|x\|_{3}\leq \dots$ holds. This sequence also determines the topology of the space. Every Fr\'echet space can be given a grading, and a graded Fr\'echet space is simply a Fr\'echet space equipped with such a grading. In this paper, we will assume that all Fr\'echet spaces discussed are graded.

A pair of graded Fr\'echet spaces $(E,F)$ is said to tame if there exists an increasing function $\sigma:\mathbb{N}\to \mathbb{N}$ such that for any continuous linear operator $T:E\to F$, there exists an $N\in \mathbb{N}$ and $C>0$ satisfying $\|Tx\|_{n}\leq C\|x\|_{\sigma(n)}$ for all $x\in E$ and $n\geq N$. A Fr\'echet space $E$ is considered tame if the pair $(E,E)$ is tame. The concept of tameness provides a way to control the continuity of operators. Dubinsky and Vogt introduced the tame Fr\'echet spaces in \cite{DV} and used it to identify a basis in complemented subspaces of certain infinite-dimensional power series spaces.

The author focused specifically on a subset of operators rather than considering all operators defined on a Fréchet space and gave the definition of the S-tameness in \cite{N} as follows:

\begin{defn} Let  $S:\mathbb{N}\to \mathbb{N}$ be a non-decreasing function. A family of linear continuous operators $\mathcal{A}\subseteq L(E,F)$ is called S-tame if for every operator $T\in \mathcal{A}$ there exist $k_{0}\in \mathbb{N}$ and  $C>0$ such that
$$\hspace{1.65in}
\|Tx\|_{k}\leq C\|x\|_{S(k)} \hspace{1in}\forall x\in E, k\geq k_{0}. $$
\end{defn}
%In \cite{N}, the author showed that S-tameness of a family can be given by considering only elements of bases similar to Lemma \ref{Crone}.
% \begin{lem}
%Let $K(a_{n,k})$ and $K(b_{n,k})$ be K\"{o}the spaces and $S:\mathbb{N}\to \mathbb{N}$ be a non-decreasing function. A family of linear continuous operators $\mathcal{A}\subseteq L(K(a_{n,k}),K(b_{n,k}))$ is S-tame if and only if for every operator $T\in \mathcal{A}$ there exist $k_{0}\in \mathbb{N}$ and $C>0$ such that
%$$\hspace{1.65in}\|Te_{n}\|_{k}\leq C
%\|e_{n}\|_{S(k)} \hspace{0.9in}\forall n\in \mathbb{N}, k\geq k_{0}.$$
%\end{lem}

We would like to note that if a family $\mathcal{A}$ of linear, continuous operators is $S_{1}$-tame and $S_{1}(n)\leq S_{2}(n)$ for sufficiently large $n\in \mathbb{N}$, then it is obvious that the family $\mathcal{A}$ is also $S_{2}$-tame. 
 
The author characterized the S-tameness of a family of operators defined by Toeplitz matrices between power series spaces in \cite{N}. Here, we will discuss the S-tameness of a family of operators defined by Hankel matrices between power series spaces.

Firstly we want to emphasize that a family of compact operators is $I$-compact where $I:\mathbb{N}\to \mathbb{N}$ is the identity. Let us assume that $\mathcal{A}$ is a family of compact operators from $K(a_{k,n})$ to $K(b_{k,n})$. Then for every $T\in A$, there exists a $m\in \mathbb{N}$ such that for all $k\geq m$ there exists a $C>0$ such that
\[\|Tx\|_{k}\leq C\|x\|_{m}\leq C\|x\|_{k}\]
for all $x\in K(a_{k,n})$. This means that the family $A$ is $I$-tame.

Now we want to address the $I$-tameness of the family of operators defined by Hankel matrices between power series spaces.
\begin{itemize}
\item[1.]  $H_{\theta}: \Lambda_{\infty}(\alpha)\to \Lambda_{\infty}(\beta)$ is compact for every $\theta\in \Lambda_{\infty}(\beta)$ by Theorem \ref{T1}. Then the family
$$\mathcal{A}=\{H_{\theta}:\Lambda_{\infty}(\alpha)\to \Lambda_{\infty}(\beta)\; | \;\theta\in \Lambda_{\infty}(\beta)\}$$
is $I$-tame.
\item[2.] $H_{\theta}: \Lambda_{1}(\alpha) \to \Lambda_{\infty}(\beta)$ is compact for every $\beta,\alpha$ satisfying the condition \ref{C1} by Theorem \ref{T2}. Then,  the family
$$\mathcal{B}=\{H_{\theta}:\Lambda_{1}(\alpha)\to \Lambda_{\infty}(\beta)\; |\; \theta\in \Lambda_{\infty}(\beta)\}$$
is $I$-tame.
\item[3.] $H_{\theta}: \Lambda_{\infty}(\alpha) \to \Lambda_{1}(\beta)$ is compact for every $\alpha$ and stable $\beta$ satisfying the condition \ref{C2} by Theorem \ref{T3}. Then, the family
$$\mathcal{C}=\{H_{\theta}:\Lambda_{\infty}(\alpha)\to \Lambda_{1}(\beta) \;|\; \theta\in \Lambda_{1}(\beta)\}$$
is $I$-tame.

\item[4.] $H_{\theta}: \Lambda_{1}(\alpha) \to \Lambda_{1}(\beta)$ is compact for every $\theta\in (\Lambda_{1}(\alpha))^{\prime}$ provided that $\Lambda_{1}(\beta)$ is a nuclear power series space of finite type by Theorem \ref{T4}. Then, the family
$$\mathcal{D}=\{H_{\theta}:\Lambda_{1}(\alpha)\to \Lambda_{1}(\beta) \;|\; \theta\in (\Lambda_{1}(\alpha))^{\prime}\}$$
is $I$-tame.
\end{itemize}

%--------------------------------------------------------------------------------------------------------------------------------------------------------------------------------------------------
\section{The Interactions of Toeplitz and Hankel Operators with Shift Operators}
In section 3, a Hankel operator $H_{\theta}:K(a_{n,k})\to K(b_{n,k})$, which we associated with a Hankel matrix
\[
\begin{pmatrix}
\theta_{0} &\theta_{1}&\theta_{2}&\theta_{3}&\cdots \\
\theta_{1}&\theta_{2}&\theta_{3}&\theta_{4}&\cdots\\
\theta_{2}&\theta_{3}&\theta_{4}&\theta_{5}&\cdots\\
\theta_{3}&\theta_{4}&\theta_{5}&\theta_{6}&\cdots\\
\vdots&\vdots&\vdots&\vdots&\ddots
\end{pmatrix}
\]
corresponding to a sequence $\theta=(\theta_{n})_{n\in \mathbb{N}_{0}}$ was defined as
\[H_{\theta}e_{n}=(\theta_{n-1},\theta_{n},\theta_{n+1},\cdots)=\sum^{\infty}_{j=1}\theta_{j+n-2}e_{j}\]
for all $n\in \mathbb{N}$. We discussed the necessary conditions for such an operator to be well-defined between power series spaces in section 3.

Similarly, in \cite{N}, a Toeplitz operator $\widehat{T}_{\theta}:K(a_{n,k})\to K(b_{n,k})$ whose associate matrix is a lower triangular Toeplitz matrix
$$
\begin{pmatrix}
\theta_{0} &0&0&0&\cdots \\
\theta_{1}&\theta_{0}&0&0&\cdots\\
\theta_{2}&\theta_{1}&\theta_{0}&0&\cdots\\
\theta_{3}&\theta_{2}&\theta_{1}&\theta_{0}&\cdots\\
\vdots&\vdots&\vdots&\vdots&\ddots
\end{pmatrix}$$
corresponding to a sequence $\theta=(\theta_{n})_{n\in \mathbb{N}_{0}}$ was defined as
$$\widehat{T}_{\theta}e_{n}=(0,\cdots,0,\theta_{0},\theta_{1},\theta_{2},\cdots)=\sum^{\infty}_{j=n}\theta_{j-n}e_{j}$$
for all $n\in \mathbb{N}$. The necessary conditions for the well-definedness, continuity, and compactness of these operators between power series spaces were given in \cite{N}.

In this section, some properties of shift operators defined between power series spaces will be discussed using the Hankel and Toeplitz operators.

The backward shift operator is defined as 
\[B:\Lambda_{r}(\alpha)\to \Lambda_{r}(\alpha), \hspace{0.2in} B(\theta)=(\theta_{n+1})_{n\in \mathbb{N}},\]
and the forward shift operator is defined as \[F:\Lambda_{r}(\alpha)\to \Lambda_{r}(\alpha), \hspace{0.2in}  F(\theta)=(\theta_{n-1})_{n\in \mathbb{N}}\]
where we assume that $\theta_{-n}=0$ for all $n\in \mathbb{N}$ and  here $r\in \{1,\infty\}$. These operators are well-defined and continuous in the case that $\alpha$ is a weakly-stable exponent sequence, that is, $\displaystyle\limsup_{n\in \mathbb{N}} \frac{\alpha_{n+1}}{\alpha_{n}}<\infty$. We recommend \cite{Ka} for more detailed information about shift operators on Köthe spaces.

We will proceed assuming that the sequence $\alpha$ is stable, but in some cases, this assumption is not necessary. We have the following relation operators $B$ and $F$ with $\widehat{T}_{\theta}$ and $H_{\theta}$
\begin{equation}\label{1} F^{n}(\theta)=(0,\cdots,0,\theta_{0},\theta_{1},\theta_{2},\cdots)=\widehat{T}_{\theta}(e_{n+1})\end{equation}
and
\begin{equation}\label{2} B^{n}(\theta)=(\theta_{n-1},\theta_{n},\theta_{n+1},\cdots)=H_{\theta}(e_{n+1})\end{equation}
for $n\in \mathbb{N}$.

%\begin{defn} For an $T\in L(E)$, a vector $x\in E$ is said to be \textbf{cyclic vector} for $T$ if the orbit
%$\text{orb}(T,x)=\{T^{n}x: n\geq 0\}$ has a dense linear span in $E$.
%\end{defn}

%\begin{itemize}
%\item Since $\text{orb}(F,\theta)=\big\{F^{n}(\theta): n\geq 0 \big \}=\big \{\widehat{T}_{\theta}(e_{n}): n\in \mathbb{N}\big \}$, we have that
%\[ \overline{\text{span}(\text{orb}(F,\theta))}=\text{range}(\widehat{T}_{\theta}).\]
%\item Since $\text{orb}(B,\theta)=\big\{B^{n}(\theta): n\geq 0 \big \}=\big \{H_{\theta}(e_{n}): n\in \mathbb{N}\big \}$,
%we have that
%\[ \overline{\text{span}(\text{orb}(B,\theta))}=\text{range}(H_{\theta}).\]
%\end{itemize}

\begin{defn}Let $T$ be a continuous linear operator on a Fr\'echet space $E$. 
The n-th Ces\'aro mean is  
\[ T^{[n]}:=\frac{1}{n}\sum^{n}_{m=1} T^{m}.\]
T is said to be \textbf{mean ergodic} if the limits $\lim_{n\to \infty} T^{[n]}x$, $x\in E$,
exists in $E$. T is said to be \textbf{Ces\'aro bounded} if the family $\{T^{[n]}:n\in \mathbb{N}\}$ is an equicontinuous subset of $L(E)$.
\end{defn}

If $E$ be a Montel Fr\'echet space, then $T$ is mean ergodic if and only if $T$ is Cesàro bounded and $\lim_{n\to\infty} \frac{1}{n}T^n x=0$ for every $x\in E$ by Theorem 2.5 of \cite{Ka}.

By using equations \ref{1} and \ref{2}, we have that  
\begin{equation*}
\begin{split}
\lim_{n\to \infty} F^{[n]}(\theta)&=\lim_{n\to \infty} \frac{1}{n}\sum^{n}_{m=1} F^{m}(\theta)= \lim_{n\to \infty} \frac{1}{n}\sum^{n}_{m=1} \widehat{T}_{\theta}(e_{m+1})\\ &=\lim_{n\to \infty} \widehat{T}_{\theta}\bigg(\sum^{n}_{m=1} \frac{1}{n}e_{m}\bigg )=\widehat{T}_{\theta}(0)=0
\end{split}
\end{equation*}
and 
\begin{equation*}
\begin{split}
\lim_{n\to \infty} B^{[n]}(\theta)&=\lim_{n\to \infty} \frac{1}{n}\sum^{n}_{m=1} B^{m}(\theta)= \lim_{n\to \infty} \frac{1}{n}\sum^{n}_{m=1} H_{\theta}(e_{m+1})\\ &=\lim_{n\to \infty} H_{\theta}\bigg(\sum^{n}_{m=1} \frac{1}{n}e_{m}\bigg )=H_{\theta}(0)=0
\end{split}
\end{equation*} 
for all $\theta\in \Lambda_{r}(\alpha)$, $r=1,\infty$.
This means that $F$ and $B$ are mean ergodic and hence $F$ and $B$ are Ces\'aro bounded for all $\theta\in \Lambda_{r}(\alpha)$, $r=1,\infty$. 

\begin{prop}
Forward shift operator F and Backward shift operator $B$ defined on  $\Lambda_{r}(\alpha)$, $r=1, \infty$, for stable $\alpha$ are 
mean ergodic and  Ces\'aro bounded. 
\end{prop}

We again recommend \cite{Ka} for the mean ergodicity of weighted shift operators on Köthe spaces.
%%%%%%%%%%%%%%%%%%%%%%%%%%%%%%%%%%%%%%%%%%%%%%%%%%%%%%%%%%%%%%%%%%%%%%%%%%%%%%%%
\subsection*{Acknowledgment}
The results in this paper were obtained while the author visited at University of Toledo. The author would like to thank TUBITAK for their support.

% ------------------------------------------------------------------------
\end{document}